\documentclass[10pt]{article}
\usepackage{amsfonts}
\usepackage{amssymb}
\usepackage{amsthm}
\usepackage{amsmath}
\usepackage{latexsym}
\usepackage{amssymb,color}

\definecolor{c20}{rgb}{0.,0.7,0.}
\definecolor{c30}{rgb}{0.,0.,1.}
\definecolor{c40}{rgb}{1,0.1,0.7}
\definecolor{c50}{rgb}{1,0,0}

\newtheorem{theorem}{Theorem}[section]
\newtheorem{lemma}{Lemma}[section]
\newtheorem{corollary}{Corollary}[section]

\numberwithin{equation}{section}
\begin{document}

\centerline{\bf \large Tail Behaviour of Weighted Sums of Order Statistics of Dependent Risks}

\bigskip
\centerline{Enkelejd Hashorva\footnote{Faculty of Business and
Economics (HEC Lausanne), University of Lausanne,  Lausanne, Switzerland}
 \  \&  Jinzhu Li \footnote{School of Mathematical Science and LPMC, 
Nankai University, Tianjin 300071, P.R. China}
}


{\bf Abstract}: 
Let $X_{1},\ldots ,X_{n}$ be $n$ real-valued dependent random variables.
With motivation from Mitra and Resnick (2009), we derive the tail asymptotic
expansion for the weighted sum of order statistics $X_{1:n}\leq \cdots \leq
X_{n:n}$ of $X_{1},\ldots ,X_{n}$ under the general case in which the
distribution function of $X_{n:n}$ is long-tailed or rapidly varying and $%
X_{1},\ldots ,X_{n}$ may not be comparable in terms of their tail
probability. We also present two examples and an application of our results
in risk theory.



\baselineskip15pt

\section{Introduction}

In numerous finance, insurance and risk management applications, a
prevailing model for the maximum and the total sum of dependent risks is the
lognormal one; see, e.g., Foss and Richards (2010), Asmussen et al.\ (2011),
Gulisashvili and Tankov (2013), Kortschak and Hashorva (2013) and Embrechts
et al.\ (2014). The asymptotic tail behaviour of the total sum (or
aggregated risk) of lognormal based models has been first derived in
Asmussen and Rojas-Nandayapa (2008). A key characteristic of lognormal risks
is that they are rapidly varying. By resorting to extreme value theory,
Mitra and Resnick (2009) offered a new methodology for the investigation of
the tail asymptotics of the total sum of rapidly varying risks. Given the
fact that in applications risks are almost always dependent, the
aforementioned paper constitutes a significant achievement in understanding
the extremal behaviour of the maximum and the total sum of dependent risks.
In particular, for dependent nonnegative random variables (risks) $%
X_{1},\ldots ,X_{n}$ such that $\lim_{x\rightarrow \infty }\mathbb{P}%
(X_{i}>x)/\mathbb{P}(X_{1}>x)=\lambda _{i}\in \lbrack 0,\infty )$ for $1\leq
i\leq n$, under some weak dependence assumptions (referred to in this paper
as Mitra-Resnick conditions), Mitra and Resnick (2009) showed that if $X_{1}$
has a distribution function in the Gumbel max-domain of attraction (see
below for the definition), then%
\begin{equation*}
\mathbb{P}(S_{n}>x)\sim \mathbb{P}(X_{n:n}>x)\sim \left(
\sum_{i=1}^{n}\lambda _{i}\right) \mathbb{P}(X_{1}>x),\qquad x\rightarrow
\infty ,
\end{equation*}%
where $S_{n}=\sum_{i=1}^{n}X_{i}$ is the aggregated risk, $X_{1:n}\leq
\cdots \leq X_{n:n}$ are the order statistics of $X_{1},\ldots ,X_{n}$, and
\textquotedblleft $\sim $\textquotedblright\ means that the ratio of the two
sides converges to $1$.

The recent contribution Asimit et al.\ (2013) showed further that, for some
positive weights $c_{0},\ldots ,c_{n-1}$, the weighted sum $%
\sum_{i=0}^{n-1}c_{i}X_{n-i:n}$ has the following asymptotic behaviour:%
\begin{equation*}
\mathbb{P}\left( \sum_{i=0}^{n-1}c_{i}X_{n-i:n}>x\right) \sim \mathbb{P}%
(c_{0}X_{n:n}>x)\sim \left( \sum_{i=1}^{n}\lambda _{i}\right) \mathbb{P}%
(c_{0}X_{1}>x),\qquad x\rightarrow \infty ,
\end{equation*}%
if the risks $X_{1},\ldots ,X_{n}$ obey the Mitra-Resnick conditions, or $%
X_{1}$ has a regularly varying distribution function and the risks are
asymptotically independent.

However, we can not obtain the tail asymptotics of $S_{n}$ or $%
\sum_{i=0}^{n-1}c_{i}X_{n-i:n}$ by the methodology given in the
aforementioned papers if there is no proportional tail-relationship among $%
X_{1},\ldots ,X_{n}$, i.e., $\lambda _{1},\ldots ,\lambda _{n}$ do not
exist. An interesting example where this is indeed the case is that of
log-normal risks with random variances presented in Example 4.1 below.

The principal goal of this contribution is to adapt the Mitra-Resnick
methodology for dependent risks which, in terms of their tail behaviour, may
not be comparable. We shall deal with both the long-tailed and rapidly
varying (real-valued) random variables allowing for three broad dependence
models. Roughly speaking, under our setup, we shall show in Theorem \ref%
{main1} below that%
\begin{equation}
\mathbb{P}\left( \sum_{i=0}^{n-1}c_{i}X_{n-i:n}>x\right) \sim \mathbb{P}%
(c_{0}X_{n:n}>x)\sim \sum_{i=1}^{n}\mathbb{P}(c_{0}X_{i}>x),\qquad
x\rightarrow \infty  \label{max-sum0}
\end{equation}%
holds uniformly for $\left( c_{0},\ldots ,c_{n-1}\right) $ in some compact
set.

The rest of this paper is organized as follows. Section 2 gives some
definitions of asymptotic theory and some preliminary results. Our main
results are presented in Section 3 followed by examples and an application
in Section 4. The proofs of all the results are relegated to Section 5.

\section{Preliminaries}

Hereafter, all limit relations hold as $x\rightarrow \infty $ unless
otherwise stated. As usual, for two positive functions $a(x)$ and $b(x)$, we
write $a(x)=o(b(x))$ or $a(x)=o(1)b(x)$ if $\lim_{x\rightarrow \infty
}a(x)/b(x)=0$. Moreover, a real-valued random variable is always assumed to
be not only concentrated on $(-\infty ,0]$. For a real-valued random
variable $X$ with distribution function $F$, we call $X$ or $F$ heavy-tailed
if the corresponding moment generating function $\widehat{F}(\delta
)=\int_{-\infty }^{\infty }\mathrm{e}^{\delta x}\mathrm{d}F(x)$ diverges to $%
\infty $ for all $\delta >0$; otherwise we call $X$ or $F$ light-tailed.

Almost all commonly used heavy-tailed distributions belong to the long-tail
class. By definition, a real-valued random variable $X$ with distribution
function $F=1-\overline{F}$ is long-tailed, denoted by $X\in \mathcal{L}$ or 
$F\in \mathcal{L}$, if $\overline{F}(x)>0$ for any $x\geq 0$ and $\overline{F%
}(x+y)\sim \overline{F}(x)$ for any $y\in \mathbb{R}$; see, e.g., Foss et
al.\ (2013). In this case, we can define two associated sets of eventually
positive functions%
\begin{equation*}
\mathcal{H}_{X}=\mathcal{H}_{F}=\{h(\cdot ):h(\cdot )\text{ satisfies
(i)--(iii)}\}
\end{equation*}%
and%
\begin{equation*}
\mathcal{H}_{X}^{\ast }=\mathcal{H}_{F}^{\ast }=\mathcal{H}_{F}\cap
\{h(\cdot ):h(x)\rightarrow \infty \},
\end{equation*}%
where properties (i)--(iii) are specified as:

\begin{itemize}
\item[\textbf{(i)}] $h(x)=o(x)$;

\item[\textbf{(ii)}] $\overline{F}(x+yh(x))\sim \overline{F}(x)$ for any $%
y\in \mathbb{R}$;

\item[\textbf{(iii)}] $h(\cdot )$ is weakly self-neglecting (introduced by
Asmussen and Foss (2014)), i.e.,%
\begin{equation*}
\limsup_{x\rightarrow \infty }\frac{h(x+yh(x))}{h(x)}<\infty ,\qquad \forall
y\in \mathbb{R}.
\end{equation*}
\end{itemize}

\noindent Property (iii) is a weakened version of the concept of
self-neglecting, which requires further that

\begin{itemize}
\item[\textbf{(iii}$^{{\protect\large \prime }}$\textbf{)}] $h\left(
x+yh(x)\right) \sim h(x)$ for any $y\in \mathbb{R}$.
\end{itemize}

A positive function $l(\cdot )$ is slowly varying if $l(xy)\sim l(x)$ for
any $y>0$. The set $\mathcal{H}_{F}^{\ast }$ is non-empty, since in view of
Lemma 4.1 of Li et al. (2010) there exists some slowly varying function $%
h(\cdot )$ (naturally self-neglecting) such that $h(x)\rightarrow \infty $
and properties (i)--(ii) hold. Additionally, $\mathcal{H}_{F}^{\ast }$ may
also contain non-slowly-varying functions. For instance, if $F$ is regularly
varying, i.e., $\overline{F}(x)\sim l(x)x^{-\alpha }\in \mathcal{R}_{-\alpha
}$ for some $\alpha \geq 0$ and some slowly varying function $l(\cdot )$,
then one can easily check that $h(x)=x^{p}(1+\beta \sin x)\in \mathcal{H}%
_{F}^{\ast }$ for any $p\in (0,1)$ and $\beta \in (-1,1)$.

A real-valued random variable $X$ with distribution function $F$ having
upper endpoint $x_{F}:=\sup \{x:F(x)<1\}$ belongs to the Gumbel max-domain
of attraction ($\mathrm{GMDA}$) if there exists some positive scaling
function $h(\cdot )$ such that%
\begin{equation}
\lim_{x\rightarrow x_{F}}\frac{\overline{F}(x+yh(x))}{\overline{F}(x)}=%
\mathrm{e}^{-y},\qquad \forall y\in \mathbb{R}.  \label{GMDA}
\end{equation}%
In this case we write $X\in \mathrm{GMDA}(h)$ or $F\in \mathrm{GMDA}(h)$.

\vskip0.3cm

\noindent \textbf{Remark 2.1. } Hereafter a scaling function $%
h(\cdot )$ of some distribution function belonging to the $\mathrm{GMDA}$
may not be the one specified in properties (i)--(iii) and (iii$^{\prime }$).
Hence, the $h(\cdot )$ in Assumption \textbf{A} is not necessarily related 
to the $h(\cdot )$ in Assumption \textbf{B} or \textbf{C} below. The unified
symbol $h(\cdot )$ for such functions is used to simplify the
writing of our assumptions and proofs below. 

\vskip0.3cm

According to extreme value theory, the normalized maxima of a random sample
with underlying distribution function in the $\mathrm{GMDA}$ converge in
distribution to a Gumbel random variable. Additionally, if $F$ belongs to
the $\mathrm{GMDA}$ with $x_{F}=\infty $, then it belongs to the class of
rapid variation specified by the relation $\lim_{x\rightarrow \infty }%
\overline{F}(xy)/\overline{F}(x)=0$ for any $y>1$. Furthermore, if $F\in 
\mathrm{GMDA}(h)$ with $x_{F}=\infty $, then $h(\cdot )$ satisfies property
(iii$^{\prime }$) mentioned above {and $h(x)=o(x)$.} See, e.g., Resnick
(1987) or Embrechts et al.\ (1997) for these well-known results.

The class of univariate distributions in the $\mathrm{GMDA}$ includes both
light-tailed and heavy-tailed distributions with exponential distributions
and heavy-tailed Weibull distributions as respective examples. On the other
hand, if $F\in \mathrm{GMDA}(h)$ and $h(x)\rightarrow \infty $ (implying $%
x_{F}=\infty $), then $F\in \mathrm{GMDA}(h)\cap \mathcal{L}$ and $%
h^{p}(\cdot )\in \mathcal{H}_{F}^{\ast }$ for any $p\in (0,1)$. Conversely, $%
F\in \mathrm{GMDA}(h)\cap \mathcal{L}$ implies $h(x)\rightarrow \infty $ by
Lemma 2.1 of Goldie and Resnick (1988). Hence, summarizing the above
analysis, we arrive at:

\begin{lemma}
\label{connection}$F\in \mathrm{GMDA}(h)$ and $h(x)\rightarrow \infty
\Longleftrightarrow F\in \mathrm{GMDA}(h)\cap \mathcal{L}\Rightarrow
h^{p}(\cdot )\in \mathcal{H}_{F}^{\ast }$ for any $p\in (0,1)$.
\end{lemma}

\section{Main Results}

Recall that $X_{1},\ldots ,X_{n}$ are $n$ real-valued dependent random
variables and $X_{1:n}\leq \cdots \leq X_{n:n}$ are the corresponding order
statistics. Enlightened by Assumptions 2.1--2.5 of Mitra and Resnick (2009),
we shall consider in this paper the following dependence structure:

\begin{itemize}
\item[$\mathbf{A.}$] $X_{n:n}\in \mathrm{GMDA}(h)$ with a distribution
function having an infinite upper endpoint. Further, it holds that%
\begin{equation}
\lim_{x\rightarrow \infty }\frac{\mathbb{P}(\left\vert X_{i}\right\vert
>th(x),X_{j}>x)}{\mathbb{P}(X_{n:n}>x)}=0\text{ for any }1\leq i\neq j\leq n%
\text{ and any }t>0,  \label{DS1}
\end{equation}%
and%
\begin{equation}
\lim_{x\rightarrow \infty }\frac{\mathbb{P}(X_{i}>Lh(x),X_{j}>Lh(x))}{%
\mathbb{P}(X_{n:n}>x)}=0\text{ for any }1\leq i<j\leq n\text{ and some }L>0.
\label{DS2}
\end{equation}
\end{itemize}

We remark that the original Mitra-Resnick conditions include relations (\ref%
{DS1}) and (\ref{DS2}) with the maximum $X_{n:n}$ replaced by $X_{1}$.
Clearly, utilizing the maximum $X_{n:n}$ instead of $X_{1}$ relaxes the
constraint of our assumption. Additionally, compared with the original
Mitra-Resnick conditions, we drop the nonnegativity of the risks and the
tail-relationships among the risks. These improvements make our Assumption 
\textbf{A} more extensive and allow us to study some flexible dependence
structures; see Examples 4.1 and 4.2 below for details. A drawback of
Assumption \textbf{A} lies in that it is not easy to show $X_{n:n}\in 
\mathrm{GMDA}(h)$ in general. To overcome this drawback, we present Lemma
4.1 below, which gives a simple condition to verify $X_{n:n}\in \mathrm{GMDA}%
(h)$.

In addition to the dependence structure controlled by Assumption \textbf{A},
we shall also investigate the asymptotic tail behaviour of weighted sums of
order statistics of dependent risks under the following long-tail case:

\begin{itemize}
\item[$\mathbf{B.}$] $X_{n:n}\in \mathcal{L}$ and there exists some $h(\cdot
)\in \mathcal{H}_{X_{n:n}}$ such that relations (\ref{DS1}) and (\ref{DS2})
hold.
\end{itemize}

A positive function $h(\cdot )$ is dominatedly varying if $%
0<\liminf_{x\rightarrow \infty }h(xy)/h(x)\leq \limsup_{x\rightarrow \infty
}h(xy)/h(x)<\infty $ for any $y>0$. For a real-valued random variable $X$
with distribution function $F$, we write $X\in \mathcal{D}$ or $F\in 
\mathcal{D}$ if $\overline{F}$ is dominatedly varying; see, e.g., Embrechts
et al. (1997) and Cai and Tang (2004) for the details. The last case
considered in this paper is as follows:

\begin{itemize}
\item[$\mathbf{C.}$] $X_{n:n}\in \mathcal{L}\cap \mathcal{D}$ and there
exists some dominatedly varying $h(\cdot )\in \mathcal{H}_{X_{n:n}}$ such
that relation (\ref{DS1}) holds.
\end{itemize}

Now, we are ready to state our main theorem, which implies the max-sum
equivalence of $X_{1},\ldots ,X_{n}$ when $c_{0}=c_{1}=\cdots =c_{n-1}=1$.
As mentioned before, compared with Corollary 2.2 of Mitra and Resnick
(2009), it contains the long-tail case and drops the nonnegativity of $%
X_{1},\ldots ,X_{n}$ and the tail-relationships among $X_{1},\ldots ,X_{n}$.

\begin{theorem}
\label{main1}Let $X_{1},\ldots ,X_{n}$ be $n$ real-valued random variables.
If one of Assumptions $\mathbf{A}$--$\mathbf{C}$ holds, then for any $%
0<a\leq b<\infty $ and $0\leq d<\infty $ relation \eqref{max-sum0} holds
uniformly for $(c_{0},c_{1},\ldots ,c_{n-1})\in \lbrack a,b]\times \lbrack
0,d]^{n-1}$.
\end{theorem}

From the proof of Theorem \ref{main1} below, we have a corresponding result
for nonnegative $X_{1},\ldots ,X_{n}$ with real-valued weights:

\begin{corollary}
\label{C1}In addition to the conditions of Theorem \ref{main1}, if $%
X_{1},\ldots ,X_{n}$ are nonnegative, then (\ref{max-sum0}) holds uniformly
for $(c_{0},c_{1},\ldots ,c_{n-1})\in \lbrack a,b]\times \lbrack -d,d]^{n-1}$%
.
\end{corollary}

Based on Theorem \ref{main1} and Corollary \ref{C1}, by conditioning on the
values of $C_{0},\ldots ,C_{n-1}$, we obtain the following corollary, in
which the assertion under Assumption $\mathbf{C}$ generalizes Theorem 1.1 of
Yang (2014).

\begin{corollary}
\label{C2}Under the conditions of Theorem \ref{main1}, let $C_{0},\ldots
,C_{n-1}$ be $n$ arbitrarily dependent random variables independent of $%
X_{1},\ldots ,X_{n}$ such that $\mathbb{P}\left( a\leq C_{0}\leq b\right) =%
\mathbb{P}\left( 0\leq C_{i}\leq d\right) =1$ for $1\leq i\leq n-1$. Then,
we have%
\begin{equation}
\mathbb{P}\left( \sum_{i=0}^{n-1}C_{i}X_{n-i:n}>x\right) \sim \mathbb{P}%
(C_{0}X_{n:n}>x)\sim \sum_{i=1}^{n}\mathbb{P}(C_{0}X_{i}>x).
\label{random weights}
\end{equation}%
If further $X_{1},\ldots ,X_{n}$ are nonnegative, then (\ref{random weights}%
) holds given that $\mathbb{P}\left( a\leq C_{0}\leq b\right) =\mathbb{P}%
\left( \left\vert C_{i}\right\vert \leq d\right) =1$ for $1\leq i\leq n-1$.
\end{corollary}

For $n$ mutually independent random variables $X_{1},\ldots ,X_{n}$, it is
of interest to seek conditions such that they are max-sum equivalent; see
Embrechts and Goldie (1980), Cai and Tang (2004), Geluk (2009), Li and Tang
(2010) and the references therein. This is connected to the well-known
principle of a single big jump in risk theory; see, e.g., Embrechts et al.
(1997) or Foss et al.\ (2007).

Next, we discuss Assumptions $\mathbf{A}$--$\mathbf{C}$ for independent $%
X_{1},\ldots ,X_{n}$. By Lemma \ref{connection} and Lemma \ref{PC}(b) below,
in the independence case our Assumptions $\mathbf{A}$--$\mathbf{C}$ have
their respective counterparts as follows:

\begin{itemize}
\item[$\mathbf{A}^{\prime }\mathbf{.}$] $X_{n:n}\in \mathrm{GMDA}(h)\cap 
\mathcal{L}$ and%
\begin{equation}
\lim_{x\rightarrow \infty }\frac{\mathbb{P}\left( X_{i}>Lh(x)\right) \mathbb{%
P}\left( X_{j}>Lh(x)\right) }{\mathbb{P}(X_{n:n}>x)}=0\text{ for any }1\leq
i\neq j\leq n\text{ and some }L>0.  \tag{3.2$^{\prime }$}  \label{IDS2}
\end{equation}

\item[$\mathbf{B}^{\prime }\mathbf{.}$] $X_{n:n}\in \mathcal{L}$ and there
exists some $h(\cdot )\in \mathcal{H}_{X_{n:n}}^{\ast }$ such that relation (%
\ref{IDS2}) holds.

\item[$\mathbf{C}^{\prime }\mathbf{.}$] $X_{n:n}\in \mathcal{L}\cap \mathcal{%
D}$.
\end{itemize}

\noindent Hence, as mentioned before, our Theorem \ref{main1} indicates the
max-sum equivalence of independent random variables meeting one of
Assumptions $\mathbf{A}^{\prime }$--$\mathbf{C}^{\prime }$. However, the
assertion under Assumption $\mathbf{C}^{\prime }$ is covered by a more
general existing result presented in Theorem 1 of Li and Tang (2010). We
conclude the assertions under Assumptions $\mathbf{A}^{\prime }$ and $%
\mathbf{B}^{\prime }$ by the following corollary.

\begin{corollary}
\label{C3}Let $X_{1},\ldots ,X_{n}$ be $n$ real-valued and mutually
independent random variables with $X_{n:n}\in \mathcal{L}$. If either (a) $%
X_{n:n}\in \mathrm{GMDA}(h)\cap \mathcal{L}$ and (\ref{IDS2}) holds or (b)
there exists some $h(\cdot )\in \mathcal{H}_{X_{n:n}}^{\ast }$ such that (%
\ref{IDS2}) holds, then $X_{1},\ldots ,X_{n}$ are max-sum equivalent, i.e.,%
\begin{equation*}
\mathbb{P}\left( \sum_{i=1}^{n}X_{i}>x\right) \sim \sum_{i=1}^{n}\mathbb{P}%
(X_{i}>x).
\end{equation*}
\end{corollary}

\section{Examples and Application}

We begin this section with a fundamental lemma, which will be applied in the
verification of the examples given below. Actually, this lemma provides a
way to verify that $X_{n:n}$ belongs to some distribution class considered
in this paper.

\begin{lemma}
\label{verify}Let $X_{1},\ldots ,X_{n}$ be $n$ real-valued random variables.
Assume that%
\begin{equation}
\lim_{x\rightarrow \infty }\frac{\mathbb{P}(X_{i}>x,X_{j}>x)}{\mathbb{P}%
(X_{n:n}>x)}=0\text{ for any }1\leq i<j\leq n.  \label{v1}
\end{equation}

(a) If $X_{i}\in \mathrm{GMDA}(h_{i})$ with $h_{i}(x)\rightarrow \infty $
and $h_{i}(x)\sim h_{1}(x)$\ for every $1\leq i\leq n$, then $X_{n:n}\in 
\mathrm{GMDA}(h_{1})$.

(b) If $X_{i}\in \mathcal{L}$ or $\mathcal{D}$ or $\mathcal{R}_{-\alpha }$
(for some $\alpha >0$) for every $1\leq i\leq n$, then $X_{n:n}\in \mathcal{L%
}$ or $\mathcal{D}$ or $\mathcal{R}_{-\alpha }$, respectively.
\end{lemma}

\proof Note the fact that%
\begin{equation}
\sum_{i=1}^{n}\mathbb{P}(X_{i}>x)-\sum_{1\leq i<j\leq n}\mathbb{P}%
(X_{i}>x,X_{j}>x)\leq \mathbb{P}(X_{n:n}>x)\leq \sum_{i=1}^{n}\mathbb{P}%
(X_{i}>x).  \label{fact}
\end{equation}%
Hence, relation (\ref{v1}) implies that%
\begin{equation}
\mathbb{P}(X_{n:n}>x)\sim \sum_{i=1}^{n}\mathbb{P}(X_{i}>x).  \label{ms}
\end{equation}%
In view of Resnick (1987), relation (\ref{GMDA}) holds locally uniformly for 
$y\in \mathbb{R}$. Hence, for assertion (a), using relation (\ref{ms}) and $%
X_{i}\in \mathrm{GMDA}(h_{i})$ with $h_{i}(x)\sim h_{1}(x)$\ for every $%
1\leq i\leq n$, it is easy to obtain relation (\ref{GMDA}) with $h(x)$
replaced by $h_{1}(x)$ for the tail probability of $X_{n:n}$. Assertion (b)
immediately follows from relation (\ref{ms}) and the definitions of the
classes $\mathcal{L}$, $\mathcal{D}$, and $\mathcal{R}_{-\alpha }$.%
\endproof

Next, we present two examples satisfying Assumption $\mathbf{A}$ and
Assumptions $\mathbf{B}$--$\mathbf{C}$, respectively. In both examples there
is no necessary tail-relationship among the random variables $X_{1},\ldots
,X_{n}$.

\vskip0.3cm

\noindent \textbf{Example 4.1. }Let $(Y_{1},\ldots ,Y_{n})$ be a
multivariate normal random vector with mean vector $(\mu _{1},\ldots ,\mu
_{n})$ and covariance matrix $\left( \rho _{ij}\sigma _{i}\sigma _{j}\right)
_{n\times n}$, where $\sigma _{i}>0,\rho _{ii}=1$ for $1\leq i\leq n$ and $%
-1<\rho _{ij}=\rho _{ji}<1$ for $1\leq i\neq j\leq n$. Let further $%
W_{1},\ldots ,W_{n}$, independent of $(Y_{1},\ldots ,Y_{n})$, be $n$
nonnegative and arbitrarily dependent random variables with finite and
positive upper endpoints. Now we verify Assumption $\mathbf{A}$ for the
random variables $X_{1}=\mathrm{e}^{W_{1}Y_{1}},\ldots ,X_{n}=\mathrm{e}%
^{W_{n}Y_{n}}$.

Since $\mu _{i}$ and $\sigma _{i}$, $1\leq i\leq n$, are arbitrarily fixed,
we simply assume that all the upper endpoints of $W_{1},\ldots ,W_{n}$ equal
to $1$ without loss of generality. We first verify that $X_{n:n}\in \mathrm{%
GMDA}(h)$ for some $h(\cdot )$. For any $1\leq i\neq j\leq n$, assume
without loss of generality that $\sigma _{i}\geq \sigma _{j}$. It is clear
for positive $x$ that%
\begin{eqnarray*}
\frac{\mathbb{P}(X_{i}>x,X_{j}>x)}{\mathbb{P}(X_{n:n}>x)} &\leq &\frac{%
\mathbb{P}(W_{i}Y_{i}>\log x,W_{j}Y_{j}>\log x)}{\mathbb{P}(W_{i}Y_{i}>\log
x)} \\
&\leq &\frac{\mathbb{P}(Y_{i}+Y_{j}>2\log x)}{\mathbb{P}(W_{i}Y_{i}>\log x)}
\\
&\leq &\left. \overline{\Phi }\left( \frac{2}{\sqrt{2\left( 1+\rho
_{ij}\right) }}\frac{\log x-\left( \mu _{i}+\mu _{j}\right) /2}{\sigma _{i}}%
\right) \right/ \mathbb{P}(W_{i}Y_{i}>\log x),
\end{eqnarray*}%
where $\Phi (\cdot )$ is the standard normal distribution function. It
follows from Lemma A.3 of Tang and Tsitsiashvili (2004) that, for any $w\in
(0,1)$,%
\begin{equation}
\mathbb{P}(W_{i}Y_{i}>x)\sim \mathbb{P}(W_{i}Y_{i}>x,W_{i}>w)\geq \mathbb{P}%
\left( Y_{i}>\frac{x}{w}\right) \mathbb{P}\left( W_{i}>w\right) =\overline{%
\Phi }\left( \frac{x-w\mu _{i}}{w\sigma _{i}}\right) \mathbb{P}\left(
W_{i}>w\right) .  \label{key}
\end{equation}%
Combining the above estimates and choosing $w>\left. \sqrt{2\left( 1+\rho
_{ij}\right) }\right/ 2$ lead to (\ref{v1}), which implies that%
\begin{equation}
\mathbb{P}(X_{n:n}>x)\sim \sum_{i=1}^{n}\mathbb{P}(X_{i}>x)\sim \sum_{i\in
\Lambda }\mathbb{P}(X_{i}>x),  \label{group}
\end{equation}%
where 
\begin{equation}
\Lambda =\left\{ i:\sigma _{i}=\sigma =\max\limits_{1\leq j\leq n}\sigma
_{j},\ \mu _{i}=\mu =\max\limits_{j:\sigma _{j}=\sigma }\mu _{j}\right\}
\label{set}
\end{equation}%
and in the last step of (\ref{group}) we used relation (\ref{key}) again. In
view of Theorem 1.1 of Hashorva and Weng (2014) (or the last sentence in the
first paragraph of their Section 2), $X_{i}\in \mathrm{GMDA}(h)$ for $i\in
\Lambda $ with the common auxiliary function $h(\cdot )$ given by%
\begin{equation*}
h(x)=\frac{\sigma ^{2}x}{\log x-\mu }.
\end{equation*}%
Hence, by Lemma \ref{verify}(a), $X_{n:n}\in \mathrm{GMDA}(h)$ with the same 
$h(\cdot )$ as above. Then, using the similar procedures as in Example 3.5
of Mitra and Resnick (2009), we can verify (\ref{DS1}) and (\ref{DS2}) with
such $h(\cdot )$. This establishes the validity of Assumption $\mathbf{A}$%
.\hfill $\Box $

For $\left( X_{1},\ldots ,X_{n}\right) $ following a multivariate lognormal
distribution as in Example 4.1 with $W_{i}\equiv 1$ for $1\leq i\leq n$,
Asmussen and Rojas-Nandayapa (2008) gave in their Theorem 1 a precise
asymptotic expansion for $\mathbb{P}\left( \sum_{i=1}^{n}X_{i}>x\right) $;
see Hashorva (2013) for some generalizations. Clearly, Example 4.1 indicates
that their result is an immediate consequence of our Theorem \ref{main1}.

\vskip0.3cm

\noindent \textbf{Example 4.2. }Consider the real-valued random variables $%
X_{1},\ldots ,X_{n}$ with distribution functions $F_{1}\in \mathcal{R}%
_{-\alpha },\ldots ,F_{n}\in \mathcal{R}_{-\alpha }$ for some $\alpha >0$.
Impose on $\left( X_{1},\ldots ,X_{n}\right) $ a multivariate
Farlie-Gumbel-Morgenstern copula (see, e.g., Hashorva and H\"{u}sler
(1999)), which implies%
\begin{equation}
\mathbb{P}(X_{1}\leq x_{1},\ldots ,X_{n}\leq
x_{n})=\prod_{i=1}^{n}F_{i}(x_{i})\left( 1+\sum_{k=2}^{n}\sum_{1\leq
j_{1}<\cdots <j_{k}\leq n}\theta _{j_{1}\cdots j_{k}}\overline{F}%
_{j_{1}}(x_{j_{1}})\cdots \overline{F}_{j_{k}}(x_{j_{k}})\right) ,
\label{FGM}
\end{equation}%
where $\left\vert \theta _{j_{1}\cdots j_{k}}\right\vert \leq 1$ are some
real numbers such that the right-hand side of (\ref{FGM}) is a proper
multivariate distribution function. We verify Assumptions $\mathbf{B}$--$%
\mathbf{C}$ for $X_{1},\ldots ,X_{n}$.

In this case, it is known that, for any $1\leq i\neq j\leq n$,%
\begin{equation}
\mathbb{P}(X_{i}>x_{i},X_{j}>x_{j})=\overline{F}_{i}(x_{i})\overline{F}%
_{j}(x_{j})\left( 1+\theta _{ij}F_{i}(x_{i})F_{j}(x_{j})\right) .  \label{2d}
\end{equation}%
Relation (\ref{2d}) obviously implies (\ref{v1}). Hence, by Lemma \ref%
{verify}(b), $X_{n:n}\in \mathcal{R}_{-\alpha }\subset \mathcal{L}\cap 
\mathcal{D}\subset \mathcal{L}$. Let $h(x)=x^{p}\in \mathcal{H}_{X_{n:n}}$
for some $p\in \left( 1/2,1\right) $. Clearly, $h(\cdot )$ is dominatedly
varying and (\ref{DS1}) holds in view of (\ref{2d}). To obtain (\ref{DS2}),
for every $1\leq i\leq n$ we write $\overline{F}_{i}(x)\sim
l_{i}(x)x^{-\alpha }$ with some slowly varying function $l_{i}(\cdot )$.
Relation (\ref{2d}) gives that, for any $1\leq i\neq j\leq n$,%
\begin{equation*}
\mathbb{P}(X_{i}>h(x),X_{j}>h(x))\sim l_{i}(x^{p})l_{j}(x^{p})x^{-2p\alpha
}\left( 1+\theta _{ij}\right) =o\left( \overline{F}_{i}(x)\right) =o(1)%
\mathbb{P}(X_{n:n}>x).
\end{equation*}%
Hence, both Assumptions $\mathbf{B}$ and $\mathbf{C}$ hold.\hfill $\Box $

Next, we present an application of our Theorem \ref{main1} in risk theory.
Let $X_{1},\ldots ,X_{n}$ be $n$ insurance risks (claims), which are
naturally nonnegative. One of popular risk measures based on the conditional
tail expectation (CTE) is defined as (recall $S_{n}=\sum_{i=1}^{n}X_{i}$)%
\begin{equation}
\mathbb{E}\left( \left. X_{i}\right\vert S_{n}>\mathrm{VaR}%
_{q}(S_{n})\right) ,\qquad 1\leq i\leq n,  \label{CTE}
\end{equation}%
where $q\in (0,1)$ and $\mathrm{VaR}_{q}(S_{n})=\inf \{x:\mathbb{P}%
(S_{n}\leq x)\geq q\}$. 

The recent contribution Asimit et al. (2011) and Zhu and Li (2012) proposed
to study the asymptotic behaviour of (\ref{CTE}) as $q\rightarrow 1$
(equivalently, $\mathrm{VaR}_{q}(S_{n})\rightarrow \infty $).

We consider a slightly broader risk measure defined by%
\begin{equation}
\mathbb{E}\left( \left. \sum_{i\in \Omega }X_{i}\right\vert S_{n}>\mathrm{VaR%
}_{q}(S_{n})\right) =\sum_{i\in \Omega }\mathbb{E}\left( \left.
X_{i}\right\vert S_{n}>\mathrm{VaR}_{q}(S_{n})\right) ,\qquad \varnothing
\neq \Omega \subset \{1,\ldots ,n\}.  \label{CTE2}
\end{equation}%
The main motivation for (\ref{CTE2}) is that risks are usually grouped and,
for risk management purposes, it is important to calculate the CTE for a
group of risks. Applying our Theorem \ref{main1}, we can obtain a pair of
asymptotic lower and upper bounds for (\ref{CTE2}) as $q\rightarrow 1$ under
Assumption $\mathbf{A}$, which extends Theorem 3.3 of Asimit et al. (2011);
see Remark 4.1 below.

\begin{theorem}
\label{main2}Let $X_{1},\ldots ,X_{n}$ be $n$ nonnegative random variables
satisfying Assumption $\mathbf{A}$. For every $\varnothing \neq \Omega
\subset \{1,\ldots ,n\}$, write%
\begin{equation*}
0\leq u=\liminf_{x\rightarrow \infty }\frac{\sum_{i\in \Omega }\mathbb{P}%
(X_{i}>x)}{\sum_{i=1}^{n}\mathbb{P}(X_{i}>x)}\leq \limsup_{x\rightarrow
\infty }\frac{\sum_{i\in \Omega }\mathbb{P}(X_{i}>x)}{\sum_{i=1}^{n}\mathbb{P%
}(X_{i}>x)}=U\leq 1.
\end{equation*}%
Then it holds that%
\begin{equation}
u\leq \liminf_{q\rightarrow 1}\frac{\mathbb{E}\left( \left. \sum_{i\in
\Omega }X_{i}\right\vert S_{n}>\mathrm{VaR}_{q}(S_{n})\right) }{\mathrm{VaR}%
_{q}(S_{n})}\leq \limsup_{q\rightarrow 1}\frac{\mathbb{E}\left( \left.
\sum_{i\in \Omega }X_{i}\right\vert S_{n}>\mathrm{VaR}_{q}(S_{n})\right) }{%
\mathrm{VaR}_{q}(S_{n})}\leq U.  \label{bounds}
\end{equation}
\end{theorem}

In view of Lemma 2.4 of Asimit et al. (2011) and Theorem \ref{main1}, we
obtain 
\begin{equation*}
\mathrm{VaR}_{q}(S_{n})\sim \mathrm{VaR}_{q}(X_{n:n}),\qquad q\rightarrow 1,
\end{equation*}%
and hence the relations in (\ref{bounds}) also hold with the denominator $%
\mathrm{VaR}_{q}(S_{n})$ replaced by $\mathrm{VaR}_{q}(X_{n:n})$.
Additionally, if we further assume that each $\lim_{x\rightarrow \infty }%
\mathbb{P}(X_{i}>x)/\mathbb{P}(X_{1}>x)$ exists for $1\leq i\leq n$ like in
Assumption 3.3 of Asimit et al. (2011), then our (\ref{bounds}) with $\Omega
\in \left\{ \{i\}:1\leq i\leq n\right\} $ reduces to their precise
asymptotic formula (3.30).

\vskip0.3cm

\noindent \textbf{Remark 4.1. }In dealing with grouped risks without
comparable tails, Theorem \ref{main2} possesses its own advantages. To see
this point, recall the generalized log-normal risks $X_{i}=\mathrm{e}%
^{W_{i}Y_{i}}$ for $1\leq i\leq n$ given in Example 4.1. In this case we are
naturally concerned with the group of dominating risks, i.e., $\left\{
X_{i}:i\in \Lambda \right\} $; see (\ref{set}). Since in general there is no
proportional\ tail-relationship among such dominating risks, Theorem 3.3 of
Asimit et al. (2011) cannot be utilized to derive the asymptotics for $%
\mathbb{E}\left( \left. \sum_{i\in \Lambda }X_{i}\right\vert S_{n}>\mathrm{%
VaR}_{q}(S_{n})\right) $. However, our Theorem \ref{main2} and relation (\ref%
{group}) give that%
\begin{equation*}
\mathbb{E}\left( \left. \sum_{i\in \Lambda }X_{i}\right\vert S_{n}>\mathrm{%
VaR}_{q}(S_{n})\right) \sim \mathrm{VaR}_{q}(S_{n}),\qquad q\rightarrow 1.
\end{equation*}

\section{Proofs}

We state first a lemma and then proceed with the proofs of our main results.

\begin{lemma}
\label{PC}Let $X_{1},\ldots ,X_{n}$ be $n$ real-valued random variables.
Assume that the distribution function of $X_{n:n}$ has an infinite upper
endpoint.

(a) Relation (\ref{DS1}) holds if and only if%
\begin{equation*}
\left[ \left. \frac{X_{k:n}}{h(x)}\right\vert (X_{n:n}>x)\right] \overset{%
\mathrm{p}}{\rightarrow }0,\qquad 1\leq k\leq n-1,
\end{equation*}%
where \textquotedblleft $\overset{\mathrm{p}}{\rightarrow }$%
\textquotedblright\ means convergence in probability as $x\rightarrow \infty 
$.

(b) If further $X_{1},\ldots ,X_{n}$ are mutually independent, then relation
(\ref{DS1}) holds if and only if $h(x)\rightarrow \infty $.
\end{lemma}

\proof(a): For the \textquotedblleft if\textquotedblright\ assertion, we use
the fact that, for any $1\leq i\neq j\leq n$ and any $t>0$,%
\begin{eqnarray*}
\frac{\mathbb{P}(\left\vert X_{i}\right\vert >th(x),X_{j}>x)}{\mathbb{P}%
(X_{n:n}>x)} &\leq &\frac{\mathbb{P}(\left\vert X_{n-1:n}\right\vert
>th(x),X_{n:n}>x)+\mathbb{P}(\left\vert X_{1:n}\right\vert >th(x),X_{n:n}>x)%
}{\mathbb{P}(X_{n:n}>x)} \\
&=&\mathbb{P}\left( \left. \frac{\left\vert X_{n-1:n}\right\vert }{h(x)}%
>t\right\vert X_{n:n}>x\right) +\mathbb{P}\left( \left. \frac{\left\vert
X_{1:n}\right\vert }{h(x)}>t\right\vert X_{n:n}>x\right) .
\end{eqnarray*}%
For the \textquotedblleft only if\textquotedblright\ assertion, we note
that, for any $1\leq k\leq n-1$ and any $t>0$,%
\begin{eqnarray*}
\frac{\mathbb{P}(\left\vert X_{k:n}\right\vert >th(x),X_{n:n}>x)}{\mathbb{P}%
(X_{n:n}>x)} &\leq &\frac{\mathbb{P}\left( \bigcup_{1\leq i\neq j\leq
n}\left( \left\vert X_{i}\right\vert >th(x),X_{j}>x\right) \right) }{\mathbb{%
P}(X_{n:n}>x)} \\
&\leq &\sum_{1\leq i\neq j\leq n}\frac{\mathbb{P}\left( \left\vert
X_{i}\right\vert >th(x),X_{j}>x\right) }{\mathbb{P}(X_{n:n}>x)}.
\end{eqnarray*}%
This completes the proof of assertion (a).

(b): Under the independence condition, it is clear that relation (\ref{ms})
holds. Thus, the \textquotedblleft if\textquotedblright\ assertion is
obvious. We shall prove the \textquotedblleft only if\textquotedblright\
assertion by contradiction. Therefore, suppose that there exists some $M>0$
and positive numbers $x_{m}\rightarrow \infty $ as $m\rightarrow \infty $
such that $h(x_{m})\leq M<\infty $ for all $m$. Since $\left\vert
X_{i}\right\vert $ does not degenerate at $0$ for $1\leq i\leq n$, we can
choose $t$ small enough such that%
\begin{equation*}
\rho _{n}=:\min_{1\leq i\leq n}\mathbb{P}\left( \left\vert X_{i}\right\vert
>tM\right) >0
\end{equation*}%
Hence, with $X_{n+1}=X_{1}$, we have%
\begin{eqnarray*}
\limsup_{x\rightarrow \infty }\sum_{i=1}^{n}\frac{\mathbb{P}\left(
\left\vert X_{i+1}\right\vert >th(x)\right) \mathbb{P}\left( X_{i}>x\right) 
}{\mathbb{P}(X_{n:n}>x)} &\geq &\limsup_{m\rightarrow \infty }\sum_{i=1}^{n}%
\frac{\mathbb{P}\left( \left\vert X_{i+1}\right\vert >th(x_{m})\right) 
\mathbb{P}\left( X_{i}>x_{m}\right) }{\mathbb{P}(X_{n:n}>x_{m})} \\
&\geq &\rho _{n}\limsup_{m\rightarrow \infty }\frac{\sum_{i=1}^{n}\mathbb{P}%
\left( X_{i}>x_{m}\right) }{\mathbb{P}(X_{n:n}>x_{m})} \\
&=&\rho _{n}>0,
\end{eqnarray*}%
which contradicts relation (\ref{DS1}).\endproof

\textbf{Proof of Theorem \ref{main1}: }Without loss of generality, we only
need to prove that, uniformly for $(c_{1},\ldots ,c_{n-1})\in \lbrack
0,d]^{n-1}$,%
\begin{equation}
\mathbb{P}\left( X_{n:n}+\sum_{i=1}^{n-1}c_{i}X_{n-i:n}>x\right) \sim 
\mathbb{P}(X_{n:n}>x)\sim \sum_{i=1}^{n}\mathbb{P}(X_{i}>x).  \label{max-sum}
\end{equation}%
The second relation in (\ref{max-sum}) is just relation (\ref{ms}), which
follows from relations (\ref{DS1}) (implying relation (\ref{v1})) and (\ref%
{fact}). Hence, the second relation in (\ref{max-sum}) holds under one of
Assumptions $\mathbf{A}$--$\mathbf{C}$.

Next, we turn to the first relation in (\ref{max-sum}). For any $t>0$ and
the function $h(\cdot )$ specified in Assumption $\mathbf{A}$ or $\mathbf{B}$
or $\mathbf{C}$, we have%
\begin{eqnarray*}
\mathbb{P}\left( X_{n:n}+\sum_{i=1}^{n-1}c_{i}X_{n-i:n}>x\right) &=&\mathbb{P%
}\left( X_{n:n}+\sum_{i=1}^{n-1}c_{i}X_{n-i:n}>x,X_{n:n}\leq x-th(x)\right)
\\
&&+\mathbb{P}\left(
X_{n:n}+\sum_{i=1}^{n-1}c_{i}X_{n-i:n}>x,X_{n:n}>x-th(x)\right) \\
&=&I_{1}(\overline{c},x)+I_{2}(\overline{c},x),
\end{eqnarray*}%
where $\overline{c}\in \lbrack 0,d]^{n-1}$ denotes the real vector $%
(c_{1},\ldots ,c_{n-1})$. Recall that under Assumption $\mathbf{A}$ or $%
\mathbf{B}$ relation (\ref{DS2}) holds. Thus, in these two cases, we
estimate $I_{1}(\overline{c},x)$ as%
\begin{equation*}
I_{1}(\overline{c},x)\leq \mathbb{P}\left(
\sum_{i=1}^{n-1}c_{i}X_{n-i:n}>th(x)\right) \leq \mathbb{P}\left(
(n-1)dX_{n-1:n}>th(x)\right) .
\end{equation*}%
By relation (\ref{DS2}), for $t>(n-1)dL$, it holds uniformly for $\overline{c%
}\in \lbrack 0,d]^{n-1}$ that%
\begin{eqnarray}
I_{1}(\overline{c},x) &\leq &\mathbb{P}\left(
(n-1)dX_{n-1:n}>th(x),(n-1)dX_{n:n}>th(x)\right)  \notag \\
&\leq &\sum_{1\leq i\neq j\leq n}\mathbb{P}\left(
X_{i}>Lh(x),X_{j}>Lh(x)\right)  \notag \\
&=&o(1)\mathbb{P}(X_{n:n}>x).  \label{I1}
\end{eqnarray}%
Under Assumption $\mathbf{C}$, we deal with $I_{1}(\overline{c},x)$ as%
\begin{eqnarray*}
I_{1}(\overline{c},x) &\leq &\mathbb{P}\left(
X_{n:n}+\sum_{i=1}^{n-1}c_{i}X_{n-i:n}>x,%
\sum_{i=1}^{n-1}c_{i}X_{n-i:n}>th(x)\right) \\
&\leq &\mathbb{P}\left( n(d+1)X_{n:n}>x,(n-1)dX_{n-1:n}>th(x)\right) \\
&\leq &\sum_{1\leq i\neq j\leq n}\mathbb{P}\left( X_{i}>\frac{x}{n(d+1)}%
,X_{j}>\frac{th(x)}{(n-1)d}\right) \\
&=&\sum_{1\leq i\neq j\leq n}\mathbb{P}\left( X_{i}>\frac{x}{n(d+1)},X_{j}>%
\frac{t}{(n-1)d}\frac{h(x)}{h\left( x/n(d+1)\right) }h\left( \frac{x}{n(d+1)}%
\right) \right) .
\end{eqnarray*}%
Recalling that $h(\cdot )$ is dominatedly varying, there exists some $\delta
>0$ such that $h(x)/h\left( x/n(d+1)\right) \geq \delta $ for large $x$.
Thus, we have%
\begin{eqnarray*}
I_{1}(\overline{c},x) &\leq &\sum_{1\leq i\neq j\leq n}\mathbb{P}\left(
X_{i}>\frac{x}{n(d+1)},X_{j}>\frac{t\delta }{(n-1)d}h\left( \frac{x}{n(d+1)}%
\right) \right) \\
&=&o(1)\mathbb{P}\left( X_{n:n}>\frac{x}{n(d+1)}\right) \\
&=&o(1)\mathbb{P}(X_{n:n}>x),
\end{eqnarray*}%
where in the second and the last steps we used (\ref{DS1}) and $F\in 
\mathcal{D}$, respectively. Hence, under one of Assumptions $\mathbf{A}$--$%
\mathbf{C}$, relation (\ref{I1}) holds uniformly for $\overline{c}\in
\lbrack 0,d]^{n-1}$.

For $I_{2}(\overline{c},x)$, we further write%
\begin{eqnarray}
I_{2}(\overline{c},x) &=&\mathbb{P}\left( \left.
X_{n:n}+\sum_{i=1}^{n-1}c_{i}X_{n-i:n}>x\right\vert X_{n:n}>x-th(x)\right) 
\mathbb{P}\left( X_{n:n}>x-th(x)\right)  \notag \\
&=&J_{1}(\overline{c},x)J_{2}(x).  \label{I21}
\end{eqnarray}%
It is clear that%
\begin{equation*}
J_{1}(\overline{c},x)=\mathbb{P}\left( \left. \frac{X_{n:n}-x}{h(x)}+\frac{%
\sum_{i=1}^{n-1}c_{i}X_{n-i:n}}{h(x)}>0\right\vert X_{n:n}>x-th(x)\right) .
\end{equation*}%
Hence, by Lemma \ref{PC}(a) and property (iii) of $h(\cdot )$, we obtain
that, uniformly for $\overline{c}\in \lbrack 0,d]^{n-1}$,%
\begin{equation*}
\left[ \!\left. \frac{\sum_{i=1}^{n-1}c_{i}X_{n-i:n}}{h(x)}\right\vert
(X_{n:n}>x-th(x))\!\right] \!\!=\!\!\left[ \!\left. \frac{%
\sum_{i=1}^{n-1}c_{i}X_{n-i:n}}{h(x-th(x))}\frac{h(x-th(x))}{h(x)}%
\right\vert (X_{n:n}>x-th(x))\!\right] \overset{\mathrm{p}}{\rightarrow }0.
\end{equation*}%
Additionally, under Assumption $\mathbf{A}$ with $h(\cdot )$ satisfying
property (iii$^{\prime }$), we can derive that%
\begin{eqnarray}
&&\left[ \left. \frac{X_{n:n}-x}{h(x)}\right\vert \left(
X_{n:n}>x-th(x)\right) \right]  \notag \\
&=&\left[ \left. \left( \frac{X_{n:n}-\left( x-th(x)\right) }{h(x-th(x))}%
\frac{h(x-th(x))}{h(x)}-t\right) \right\vert \left( X_{n:n}>x-th(x)\right) %
\right] \overset{\mathrm{d}}{\rightarrow }Y-t,  \label{Y}
\end{eqnarray}%
where \textquotedblleft $\overset{\mathrm{d}}{\rightarrow }$%
\textquotedblright\ means convergence in distribution as $x\rightarrow
\infty $ and $Y$ is an exponential random variable with expectation $1$.
Note further that relation (\ref{Y}) holds with $Y=\infty $ under Assumption 
$\mathbf{B}$ or $\mathbf{C}$, because of $X_{n:n}\in \mathcal{L}$ and
property (ii). Hence, it holds uniformly for $\overline{c}\in \lbrack
0,d]^{n-1}$ that%
\begin{equation*}
\left[ \left. \left( \frac{X_{n:n}-x}{h(x)}+\frac{%
\sum_{i=1}^{n-1}c_{i}X_{n-i:n}}{h(x)}\right) \right\vert \left(
X_{n:n}>x-th(x)\right) \right] \overset{\mathrm{d}}{\rightarrow }Y-t,
\end{equation*}%
which implies that, uniformly for $\overline{c}\in \lbrack 0,d]^{n-1}$,%
\begin{equation}
\lim_{x\rightarrow \infty }J_{1}(\overline{c},x)=\mathbb{P}\left(
Y-t>0\right) =\left\{ 
\begin{array}{l}
\mathrm{e}^{-t},\qquad \text{under Assumption }\mathbf{A} \\ 
\\ 
1,\qquad \quad \text{under Assumption }\mathbf{B}\text{ or }\mathbf{C}%
\end{array}%
\right. .  \label{I22}
\end{equation}%
On the other hand, we have%
\begin{equation}
\lim_{x\rightarrow \infty }\frac{J_{2}(x)}{\mathbb{P}(X_{n:n}>x)}=\left\{ 
\begin{array}{l}
\mathrm{e}^{t},\quad \,\,\,\,\text{under Assumption }\mathbf{A} \\ 
\\ 
1,\qquad \text{under Assumption }\mathbf{B}\text{ or }\mathbf{C}%
\end{array}%
\right. .  \label{I23}
\end{equation}%
Plugging (\ref{I22}) and (\ref{I23}) into (\ref{I21}) leads to that the
relation%
\begin{equation*}
I_{2}(\overline{c},x)\sim \mathbb{P}(X_{n:n}>x)
\end{equation*}%
holds uniformly for $\overline{c}\in \lbrack 0,d]^{n-1}$ under one of
Assumptions $\mathbf{A}$--$\mathbf{C}$. This, together with (\ref{I1}),
completes the proof.\hfill $\Box $

For the proof of Theorem \ref{main2}, we shall need a crucial property of
distribution functions in the $\mathrm{GMDA}$ referred to as the
Davis-Resnick tail property; see Proposition 1.1 of Davis and Resnick (1988)
or relation (5) of Balakrishnan and Hashorva (2013). Namely, for any
distribution function $F\in \mathrm{GMDA}(h)$ with an infinite upper
endpoint, the following bound holds for all large $x$ and $\varepsilon >0$:%
\begin{equation}
\frac{\overline{F}(x+h(x)s)}{\overline{F}(x)}\leq \left( 1+\varepsilon
\right) \left( 1+\varepsilon s\right) ^{-1/\varepsilon },\qquad \forall
s\geq 0,  \label{DR}
\end{equation}%
where $h(\cdot )$ is a particular scaling function such that for all large $%
x>x_{0}$ we have $\overline{F}(x)\!=\!c(x)\!\exp \!\left(
\!-\!\int_{x_{0}}^{x}h(t)\mathrm{d}t\!\right) $ with $c(\cdot )$ a
measurable function satisfying $\lim_{x\rightarrow \infty }c(x)=c>0$. Note
in passing that any other scaling function $h^{\ast }$ such that $F\in 
\mathrm{GMDA}(h^{\ast })$ is asymptotically equivalent to $h$.

\vskip0.3cm

\textbf{Proof of Theorem \ref{main2}: }Without loss of generality, we only
prove the case of $\Omega =\left\{ 1,\ldots ,m\right\} $ for some $1\leq
m\leq n$. Following the proof of Theorem 3.3 of Asimit et al. (2011), we
write%
\begin{eqnarray*}
\mathbb{E}\left( \left. S_{m}\right\vert S_{n}>x\right) &=&\left(
\int_{0}^{h(x)}+\int_{h(x)}^{x}+\int_{x}^{\infty }\right) \mathbb{P}(\left.
S_{m}>z\right\vert S_{n}>x)\mathrm{d}z \\
&=&I_{1}(x)+I_{2}(x)+I_{3}(x).
\end{eqnarray*}%
It is clear that $I_{1}(x)\leq h(x)=o(x)$. By the change of variable $%
z=x+h(x)s$, we have%
\begin{eqnarray*}
I_{3}(x) &=&h(x)\int_{0}^{\infty }\mathbb{P}(\left. S_{m}>x+h(x)s\right\vert
S_{n}>x)\mathrm{d}s \\
&\leq &h(x)\int_{0}^{\infty }\frac{\mathbb{P}(S_{n}>x+h(x)s)}{\mathbb{P}%
(S_{n}>x)}\mathrm{d}s \\
&\sim &h(x)\int_{0}^{\infty }\frac{\mathbb{P}(X_{n:n}>x+h(x)s)}{\mathbb{P}%
(X_{n:n}>x)}\mathrm{d}s,
\end{eqnarray*}%
where in the last step we used Theorem \ref{main1}. Hence, by the
aforementioned Davis-Resnick tail property (\ref{DR}) for $X_{n:n}\in 
\mathrm{GMDA}(h)$ and the Dominated Convergence Theorem, we have%
\begin{equation*}
I_{3}(x)\leq (1+o(1))h(x)\int_{0}^{\infty }\mathrm{e}^{-s}\mathrm{d}s=o(x).
\end{equation*}%
It remains to verify that%
\begin{equation*}
\liminf_{x\rightarrow \infty }\frac{\sum_{i=1}^{m}\mathbb{P}(X_{i}>x)}{%
\sum_{i=1}^{n}\mathbb{P}(X_{i}>x)}\leq \liminf_{x\rightarrow \infty }\frac{%
I_{2}(x)}{x}\leq \limsup_{x\rightarrow \infty }\frac{I_{2}(x)}{x}\leq
\limsup_{x\rightarrow \infty }\frac{\sum_{i=1}^{m}\mathbb{P}(X_{i}>x)}{%
\sum_{i=1}^{n}\mathbb{P}(X_{i}>x)}.
\end{equation*}%
To this purpose, we further write%
\begin{eqnarray}
\liminf_{x\rightarrow \infty }\frac{I_{2}(x)}{x} &\geq
&\liminf_{x\rightarrow \infty }\frac{(x-h(x))\mathbb{P}(\left.
S_{m}>x\right\vert S_{n}>x)}{x}  \notag \\
&=&\liminf_{x\rightarrow \infty }\frac{\mathbb{P}(S_{m}>x)}{\mathbb{P}%
(S_{n}>x)}-\lim_{x\rightarrow \infty }\frac{h(x)}{x}\mathbb{P}(\left.
S_{m}>x\right\vert S_{n}>x)  \notag \\
&\geq &\liminf_{x\rightarrow \infty }\frac{\sum_{i=1}^{m}\mathbb{P}%
(X_{i}>x)-\sum_{1\leq i<j\leq m}\mathbb{P}(X_{i}>x,X_{j}>x)}{\mathbb{P}%
(S_{n}>x)}  \notag \\
&=&\liminf_{x\rightarrow \infty }\frac{\sum_{i=1}^{m}\mathbb{P}(X_{i}>x)}{%
\sum_{i=1}^{n}\mathbb{P}(X_{i}>x)},  \label{one hand}
\end{eqnarray}%
where in the last step we used Theorem \ref{main1} and Assumption $\mathbf{A}
$. Additionally, it holds that%
\begin{eqnarray*}
&&\mathbb{P}\left( S_{m}>h(x),S_{n}>x\right) \\
&=&\mathbb{P}(S_{n}>x)-\mathbb{P}(S_{m}\leq h(x),S_{n}>x) \\
&\leq &\mathbb{P}(S_{n}>x)-\mathbb{P}\left( S_{m}\leq
h(x),\bigcup\limits_{j=m+1}^{n}(X_{j}>x)\right) \\
&=&\mathbb{P}(S_{n}>x)-\mathbb{P}\left(
\bigcup\limits_{j=m+1}^{n}(X_{j}>x)\right) +\mathbb{P}\left(
S_{m}>h(x),\bigcup\limits_{j=m+1}^{n}(X_{j}>x)\right) \\
&\leq &\mathbb{P}(S_{n}>x)-\!\!\!\sum\limits_{j=m+1}^{n}\!\!\!\mathbb{P}%
\!\left( X_{j}>x\right) +\!\!\!\sum\limits_{m+1\leq j<k\leq n}\!\!\!\!%
\mathbb{P}\!\left( X_{j}>x,X_{k}>x\right) +\!\!\!\sum\limits_{1\leq i\leq
m<j\leq n}\!\!\!\!\mathbb{P}\!\left( X_{i}>\frac{h(x)}{m},X_{j}>x\right) .
\end{eqnarray*}%
It follows from Theorem \ref{main1} and Assumption $\mathbf{A}$ that the
last two terms in the above relation are asymptotically negligible compared
with $\mathbb{P}\left( S_{n}>x\right) $. Hence, by Theorem \ref{main1}
again, we have%
\begin{eqnarray*}
\limsup_{x\rightarrow \infty }\frac{I_{2}(x)}{x} &\leq
&\limsup_{x\rightarrow \infty }\frac{\mathbb{P}\left(
S_{m}>h(x),S_{n}>x\right) }{\mathbb{P}(S_{n}>x)}  \notag \\
&\leq &\limsup_{x\rightarrow \infty }\frac{\mathbb{P}(S_{n}>x)-%
\sum_{j=m+1}^{n}\mathbb{P}\left( X_{j}>x\right) }{\mathbb{P}(S_{n}>x)} 
\notag \\
&=&\limsup_{x\rightarrow \infty }\frac{\sum_{i=1}^{m}\mathbb{P}(X_{i}>x)}{%
\sum_{i=1}^{n}\mathbb{P}(X_{i}>x)}, 
\end{eqnarray*}%
which together with (\ref{one hand}) completes the proof.\hfill $\Box $

\vskip0.5cm

\noindent \textbf{Acknowledgment.} We are very much in debt to the two
reviewers and an Editor who carefully read our manuscript and
showed us how to significantly improve it. The authors were partially
supported by the Swiss National Science Foundation Project 200021-140633/1
and the project RARE -318984 (an FP7 Marie Curie IRSES Fellowship). The
second author also acknowledges the support from the National Natural
Science Foundation of China (Grant No.: 11201245).

\end{document}